\documentclass[12pt]{article}

\usepackage{amssymb,amsmath,latexsym}
\usepackage[dvips]{graphics}
\usepackage{mathrsfs}
\usepackage{color}
\usepackage{comment}

\newtheorem{thm}{Theorem}[section]
\newtheorem{prop}[thm]{Proposition}

\newtheorem{remark}{Remark}[section]

\newtheorem{cor}[thm]{Corollary}
\newtheorem{example}[thm]{Example}
\newtheorem{lemma}[thm]{Lemma}
\numberwithin{equation}{section}

\def\la{{\langle}} \def\ra{{\rangle}}
  \def\sC {{\mathcal C}}

  \def\bI {{\mathbb I}}
  
 \def\bN {{\mathbb N}} 
\def\bP {{\mathbb P}}  \def\bR {{\mathbb R}}

\def\EE{\mathbb{E}}

\def\VV{\mathscr{V}}

\def\<{\langle}
\def\>{\rangle}
\def\pf{\noindent{\bf Proof.} }

\def\VV{\mathbb{V}}
\def\P{{\bP}}

\def\qed{{\hfill $\Box$\medskip}}

\def\epsilon{\varepsilon}

\allowdisplaybreaks[4]

\def\be{\begin{equation}}

\def\ee{\end{equation}}

\def\al{\alpha}

\begin{document}
\title{\bf Pointwise upper  estimates  for  transition probabilities of continuous time random walks on graphs}
\author{   Xinxing Chen \footnote{Department of Mathematics, Shanghai Jiaotong University, Shanghai 200240,
China, chenxinx@sjtu.edu.cn.}} \maketitle
\begin{abstract}
Let $X$ be a continuous time random walk on a weighted graph. Given
the on-diagonal upper bounds of transition probabilities at two
vertices $x_1$ and $x_2$, we obtain Gaussian upper estimates for the
off-diagonal transition probability $\P_{x_1}(X_t=x_2)$ in terms of
 an adapted metric introduced  by Davies.
\end{abstract} \noindent{\bf 2000 MR subject
classification:} 60G50, 30K08

\section{Introduction}
Let $\Gamma=(\VV,\EE)$ be a connected, locally finite  graph without
double edges. The graph  $\Gamma$ can be either finite or infinite.
Let $\mu$ be an edge weight function on $\EE$, such that
$\mu_{xy}=\mu_{yx}> 0$ for each $(x,y)\in \EE$, while $\mu_{xy}=0$
for each $(x,y)\not\in \EE$. { Let $\nu_x>0$ for  $x\in \VV$.}
Denote by $X=\{X_t: t\ge 0\}$ a continuous time random walk on
$\Gamma$ with generator 
$$ \mathscr{L}f(x)=\frac{1}{\nu_x}\sum_{y\in \VV}
(f(y)-f(x))\mu_{xy}. $$ Write $\P_x$ for the probability measure of
$X$ starting from $x$.

If $\nu_x=\sum \mu_{xy}$ for all $x$, then the process $X$
  is called the {\em constant speed  random
walk} or CSRW on $\VV$. It is  a process that waits an exponential
time mean 1  at each vertex and then jumps to one of its neighbours.
If $\nu_x\equiv1$,  then  the expected waiting time of each jump may
vary. Moreover, such a process may explode in finite time.

In this paper, we fix  vertices $x_1,x_2\in \VV$ and  functions
$f_1,f_2$ on $\bR_+$ such that for any $i=1,2$ and $t\ge 0$,
\begin{equation}\label{e:diag}
\P_{x_i}(X_t=x_i)\le \frac{1}{f_i(t)}.
\end{equation}
 Our interest is, under what circumstances
$\P_{x_1}(X_t=x_2)$  will have Gaussian upper bounds.  Let
$d_{\nu}(\cdot,\cdot)$ be a metric of $\Gamma$  such that
\begin{equation}\label{e:nud}
\begin{cases}
\frac{1}{\nu_x}\sum_{y}{d_\nu}(x,y)^2\mu_{xy}\le 1  & \hbox{for
all~} x\in \VV,\\
{d_\nu}(x,y)\le 1  & \hbox{whenever~}x,y\in\VV {\rm ~and~}x\sim y.
\end{cases}
\end{equation}
Metrics satisfying (\ref{e:nud}) are called adapted metrics.  Such
metrics were introduced by Davies \cite{Da} and \cite{Da2}, and are
closely related to the intrinsic metric associated with a given
Dirichlet form. {(One might expect that  analogues  of diffusion
processes on manifolds hold using the intrinsic metrics for  random
walks on graphs, see \cite{FLW},\cite{Fl2} and \cite{SK}.)} { Fix
$A\ge 1$ and $\gamma>1$.}
 Let
$f:\bR_+\rightarrow \bR_+$.  We say that $f$ is $(A,\gamma)-$regular
on $[a,b)$, if the function $f$ is non-decreasing on $\bR_+$ and
satisfies that
\begin{equation}\label{e:regular}
\frac{f(\gamma s)}{f(s)}\le A\frac{f(\gamma t)}{f(t)}~~~{\rm for
~all~}a\le s<t<\gamma^{-1}b.
\end{equation}
In particular, if $a=0$ and $b=\infty$ then we say that $f$ is
$(A,\gamma)-$regular, which was introduced by Grigor'yan \cite{Gr}.
\begin{thm}\label{t:main4} Let $\delta\ge 1$. If each $f_i$ is $(A,\gamma)-$regular and
satisfies
   \begin{equation}\label{e:speed}
f_i(t)\le A e^{\delta t} {\rm~for~all~}t\in \bR_+,
\end{equation}
then there exist  constants $C_1$, $\theta>0$ which are independent
of $A,\gamma$ and $\delta$, such that
\begin{align*}
\P_{x_1}(X_{t}=x_2)\le \frac{C_1A^\beta
(\nu_{x_2}/\nu_{x_1})^{1/2}}{\sqrt{f_1(\alpha t)f_2(\alpha
t)}}\exp\left(-\theta\frac{d_\nu(x_1,x_2)^2}{ t}\right)~{\rm~for~
~}t\ge d_\nu(x_1,x_2),
\end{align*}
where $\alpha=\min\{(2\gamma)^{-1}, (64\delta)^{-1}\}$ and
$\beta=\lceil\frac{\log \gamma}{\log 2}\rceil$.\\
\end{thm}

The problem of getting a Gaussian upper bound from two point
estimates was introduced in the manifold case by Grigor'yan
\cite{Gr}. In  subsequent researches,  Coulhon, Grigor'yan \& Zucca
\cite{CGZ} studied the problem for discrete time random walks on
graphs, while
 Folz \cite{FL} studied it for the  continuous time random walks.
The current paper considers the same
 problem, however, it improves the result of \cite{FL} by no longer
requiring a lower bound on $\nu_x$.   The improvement comes from
imposing conditions on the transition probabilities $\P_x(X_t=x)$
instead of the heat kernels $p_t(x,x)$. Note that the transition
probabilities are invariant under the transformation from
$(\mu,\nu)$ to $(c\mu,c\nu)$, where $(c\mu)_{xy}=c\mu_{xy}$ and
$(c\nu)_{x}=c\nu_{x}$.

\begin{remark}\label{r:1612} The condition (\ref{e:speed}) is quite natural.
 Note that $ \P_x(X_t=x)\ge
\exp\left(-\frac{\mu_x}{\nu_x} t\right),$ where
$\mu_x=\sum_{y}\mu_{xy}$.  It implies that (\ref{e:speed}) holds if
$A=1$ and
$\delta=\max\{\frac{\mu_{x_1}}{\nu_{x_1}},\frac{\mu_{x_2}}{\nu_{x_2}}\}$.
In particular,  for CSRW one can  take $\delta=1$.
\end{remark}
\begin{remark} One can  also trace the values of $C_1$ and $\theta$.
Indeed, we select  $\theta=10^{-7}$ in our proof.
\end{remark}
\begin{thm}\label{t:main5} Let  $\delta\ge 1$. If each $f_i$ is $(A,\gamma)-$regular on $[T_1,T_2)$ and
satisfies
\begin{equation}\label{e:11181}
f_i(t)\le A e^{\delta t} {\rm~for~all~}t\in [T_1,T_2),
\end{equation}
then there exist  constants $C_1$, $\theta>0$ which are independent
of $A,\gamma$ and $\delta$, such that
\begin{equation}\label{e:1148}
\P_{x_1}(X_{t}=x_2)\le \frac{C_1A^\beta
(\nu_{x_2}/\nu_{x_1})^{1/2}}{\sqrt{f_1(\alpha t)f_2(\alpha
t)}}\exp\left(-\theta\frac{d_\nu(x_1,x_2)^2}{ t}\right)~~{\rm for~
}t\in [\widetilde{T_1},~T_2),
\end{equation}
where $\alpha=\min\{(2\gamma)^{-1}, (64\delta)^{-1}\}$,
$\beta=\lceil\frac{\log \gamma}{\log 2}\rceil$ and
$\widetilde{T_1}=(8\alpha^{-2} T_1^2)\vee d_{\nu}(x_1,x_2)$.
\end{thm}

\begin{remark} (1) If the growth rate of  $f_i$ is
either sub-exponential  or polynomial, then  the lower bound of
$\widetilde{T}_1$ will be improved,
see Theorems \ref{t:main2} and  \ref{t:main3}.\\
(2)Theorems \ref{t:main4} and \ref{t:main5} are potentially very
useful for random walks with random conductances. For example,
Mathieu and Remy \cite{MR} considered the CSRW on the infinite
cluster $\sC_\infty(\omega)$ and showed that $\P_x(X_t=x)\le
ct^{-d/2}$ for $t\ge N_x(\omega)$ and $x\in \sC_\infty(\omega)$.
Using Theorem \ref{t:main5}  immediately gives
$$
\P_x(X_t=y)\le
c_1t^{-d/2}\exp\left(-c_2\frac{d(x,y)^2}{t}\right),~~~~t\ge
S_{xy}(\omega)\vee d(x,y),
$$
where $S_{xy}(\omega)=64^3(N_x(\omega)^2\vee N_y(\omega)^2)$ and
$d(x,y)$ is the graph distance. A more delicate result by a
different method was obtained in \cite{MB1}.

 See 
  Balow and Chen \cite{BC} for
the new application in a deterministic graph where  volume doubling
and Poincar$\acute{e}$ inequality hold  for all sufficiently large
balls.
\end{remark}

In Section 2, we show  the Integral Maximum Principle  for a
positive subsolution function on $\bR_+\times \VV$. From this, we
get the initial  estimates  of the transition probabilities, the
case $t\le {d_\nu}(x,y)$ included. In Section 3,  we  update the
results of the previous section, under the assumption that a certain
regularity condition holds. In Section 4, we give the proof of
Theorem \ref{t:main4}. In the final section, we consider  functions
which are regular only on an interval and  have different rates of
growth; in doing so, we obtain Theorem \ref{t:main5}.

\section{Integral maximum principle}
For any functions $f,g$ on $\VV$, define
$$
\langle f, g\rangle=\sum_{x\in \VV}f(x)g(x)\nu_x.
$$
Then   $\langle\cdot,\cdot\rangle$ induces an inner product space.
Denote by  $\|\cdot\|$  the induced  norm. Let $\bI$ be an interval
of $\bR_+$.  We say that $u:\bI\times \VV\mapsto \bR_+$ is a
positive subsolution of the heat equation on $\bI\times\VV$ if
$$
\tfrac{\partial}{\partial t}u\le \mathscr{L}u
~~~\hbox{on}~~\bI\times \VV.
$$
Furthermore, we define a set of functions:
$$
\mathcal{H}(\bI)=\left\{u: u \hbox{~is a positive subsolution on~}
\bI\times\VV {\rm~ and~}
 |\{z:u(t,z)\not=0 \text{~for some~} t\in \bI\}|<\infty \right\}.
$$
  Let
$o\in B\subseteq \VV$ with $|B|<\infty$. Set
\begin{equation}\label{e:exiting}u_B(t,z)=\frac{\nu_o^{1/2}}{\nu_z}\P_o(X_t=z,~ \inf\{s\ge 0:X_s\not\in B\}>t).
\end{equation}
Then $u_B=0$ on $\bR_+\times(\VV\setminus B)$.
 Since
$\Gamma$ is a locally finite graph,
 $u_B$ is a positive
subsolution  on $\bR_+\times\VV$ and so $u_B\in \mathcal{H}(\bR_+)$.
Now we show the Integral Maximum Principle.
\begin{thm}\label{l:73.1} Let $ h\hbox{~ be  a   positive function on~}\bI\times \VV$ and  $u\in \mathcal{H}(\bI)$.   If
   for each $t\in \bI$ one has
\begin{equation}\label{e:80.0}
\frac{1}{\nu_y}\sum_x\frac{|h(t,x)-h(t,y)|^2}{4h(t,x)h(t,y)}\mu_{xy}\le
-\frac{\partial}{\partial t}\log h(t,y) \hbox{~~~for all~ } y\in
\VV,
\end{equation}
 then $J(t)=\langle u^2(t,\cdot), h(t,\cdot)\rangle$ is
non-increasing on $\bI$.
\end{thm}
\pf For brevity, we  omit the notation $t$. Set
$\nabla_{xy}g=g(t,y)-g(t,x)$  for any function $g$ on $\bI\times
\VV$ and get
\begin{align*}
&~~~\la 2u\mathscr{L}u,h  \ra=2\la uh,\mathscr{L}u\ra\\
&= -\sum_{x,y} \nabla_{xy}(uh) \cdot \nabla_{xy}u \cdot \mu_{xy}
\hbox{~~~~~~~~since~} |\{z:u(t,z)\not=0 \text{~for some~} t\in \bI\}|<\infty \\
&=-\sum_{x,y} (h(x)\nabla_{xy}u+u(y)\nabla_{xy}h) \cdot
\nabla_{xy}u \cdot \mu_{xy}\\
&=-\sum_{x,y} \big((\nabla_{xy}u)^2 h(x)+ u(y)\nabla_{xy}u \cdot
\nabla_{xy}h\big) \mu_{xy}\\
&=\sum_{x,y}\left[- \left( \sqrt{h(x)}\nabla_{xy}u+\frac{ u(y)
\nabla_{xy}h}{2\sqrt{h(x)}} \right)^2+\frac{(u(y)\nabla_{xy}h)^2}{4h(x)}\right]\mu_{xy} \hbox{~~since~} h \hbox{~is  positive} \\
&\le \sum_{x,y}u(y)^2\frac{|\nabla_{xy}h|^2}{4h(x)}\mu_{xy}\\
 &=\sum_{y}u(y)^2 \left( \sum_x\frac{|\nabla_{xy}h|^2}{4h(x)} \mu_{xy}\right) .
\end{align*}
By (\ref{e:80.0}), $\sum_x\frac{|\nabla_{xy}h|^2}{4h(x)}\mu_{xy}\le
-\nu_y\frac{\partial}{\partial t}h(y)$ and hence
$$
\la 2u\mathscr{L}u,h \ra\le - \sum_{y}u(y)^2
\nu_y\frac{\partial}{\partial t}h(y) =  -\la u^2,
\frac{\partial}{\partial t}h\ra.
$$

On the other hand, by the condition that $u$ is a positive
subsolution on $\bI\times\VV$, we have
\begin{align*}
\frac{d}{dt}J=&\frac{\partial}{\partial t}\langle
u^2,h\rangle=\langle 2u\frac{\partial}{\partial t}u,h\rangle +\la
u^2,\frac{\partial}{\partial t}h \ra \le \la 2u\mathscr{L} u,h\ra
+\la u^2, \frac{\partial}{\partial t}h\ra\le 0.
\end{align*}
Therefore,  $J$ is
non-increasing. \qed\\

Since  the metric $d_\nu$  satisfies (\ref{e:nud}), Theorem
\ref{l:73.1}  immediately implies Corollary \ref{c:2.2*}  as
follows. Define a set of functions:
\begin{align*}
 \mathcal{F}(\bI)=&\{h: h\hbox{~ is  a   positive function on~}\bI\times
 \VV \hbox{~and for each~} t\in \bI, x,y\in \VV {\rm~with~} x\sim y, \\
&~~~~~~~~~~~~~~~~\frac{|h(t,x)-h(t,y)|^2}{4h(t,x)h(t,y)}\le
-{d_\nu}(x,y)^2\frac{\partial}{\partial t}\log h(t,y)
 \Big\}.
 \end{align*}
\begin{cor}\label{c:2.2*} Let $u\in \mathcal{H}(\bI)$ and $h\in \mathcal{F}(\bI)$.
Then $J(t)=\langle u^2(t,\cdot), h(t,\cdot)\rangle$ is
non-increasing on $\bI$.
\end{cor}

Next, some useful  functions in $\mathcal{F}(\bI)$ will be given
below. Let $\rho(\cdot)$ be any nonnegative function on $\VV$ such
that
\begin{equation}\label{e:rholimitation1}|\rho(x)-\rho(y)|\le {d_\nu}(x,y)  \hbox{~~for any~} x,y\in
\VV\hbox{~with~}x\sim y.\end{equation} (In practice, one often
chooses $\rho(\cdot)=d_\nu(o, \cdot)\wedge R$ for some $o\in \VV$
and $R\ge 0$.)
\begin{lemma}\label{l:out} Let $\tau>0$.
For each $t\ge 0$ and $z\in \VV$, set
$$h(t,z)=\exp\left\{\big(\rho(z)-4^{-1}e~(t+\tau)\big)~\log\left(1\vee\frac{\rho(z)}{4^{-1}e~(t+\tau)}\right)
-\frac{t}{\tau}\right\}.$$ Then $h(t,z)\in \mathcal{F(\bR_+)}$.
\end{lemma}
\pf   We first show that for any $x\in [0,\infty)$ and
${\epsilon}\in[0,1]$,
\begin{align}\label{e:lotay}
e^{{\epsilon} x}+e^{-{\epsilon} x}-2\le&
{\epsilon}^2(e^x+e^{-x}-2)\end{align} and
\begin{align}
\label{e:lotay22} 1-e^{-\epsilon x}\ge& \epsilon(1-e^{-x}).
\end{align}
 By the Mean Value Theorem,
\begin{align*}
\frac{e^{{\epsilon} x}+e^{-{\epsilon}
x}-2}{{\epsilon}^2(e^{x}+e^{-x}-2)
}=&\frac{e^{{\epsilon}x_1}-e^{-{\epsilon}x_1}}{{\epsilon}(e^{x_1}-e^{-x_1})}=\frac{e^{{\epsilon}x_2}+e^{-{\epsilon}x_2}}{e^{x_2}+e^{-x_2}}\le
1,
\end{align*}
where $x> x_1> x_2> 0$. Consequently, (\ref{e:lotay}) holds. In the
same way,
 we can obtain (\ref{e:lotay22}).\\

 Fix
$y\sim z$ and  $\epsilon=d_\nu(y,z)$. Then $|\rho(y)-\rho(z)|\le
\epsilon\le 1$ by (\ref{e:nud}) and (\ref{e:rholimitation1}).  Write
$t^+= t+\tau$ and
$$
b=\left|(\rho(y)-4^{-1}et^+)\log\left(1\vee\frac{\rho(y)}{4^{-1}et^+}\right)
-(\rho(z)-4^{-1}et^+)\log\left(1\vee\frac{\rho(z)}{4^{-1}et^+}\right)\right|.
$$
Then
$$
\frac{|h(t,z)-h(t,y)|^2}{4h(t,z)h(t,y)} =\frac{e^b+e^{-b}-2}{4}.
$$
  We shall consider three cases.\\

  Case I:  $\rho(z),\rho(y)\le 4^{-1}et^+$. Then $b=0$ and $$
\frac{|h(t,z)-h(t,y)|^2}{4h(t,z)h(t,y)}=\frac{e^b+e^{-b}-2}{4}=0.$$

Case II:  $\rho(z),\rho(y)\ge 4^{-1}et^+$. By the Mean Value
Theorem,
\begin{align*}
b
=&|\rho(y)-\rho(z)|\left(\log\left(\frac{\xi}{4^{-1}et^+}\right)+\frac{\xi-4^{-1}et^+}{\xi}
\right),
\end{align*}
where $\xi$ is some value between $\rho(y)$ and $\rho(z)$.
Furthermore, we have $4^{-1}et^+\le \xi\le \rho(y)+\epsilon$ and
\begin{align*}
b\le& \epsilon\log\left(\frac{4\xi}{t^+}e^{-4^{-1}et^+/\xi}\right)\\
\le& \epsilon\log\left(\frac{4\xi}{t^+}\big(~1-(1-e^{-1})4^{-1}et^+/\xi~\big)\right)  ~~~~~~~~~~ \hbox{by}~(\ref{e:lotay22})  \\  
=&\epsilon\log\left(\frac{4\xi}{t^+}-e+1\right)\le
\epsilon\log\left(4\frac{\rho(y)+\epsilon}{t^+}-e+1\right).
\end{align*} As a result,
\begin{align*}
e^b+e^{-b}-2\le& \exp\left(\epsilon
\log\left(4\frac{\rho(y)+\epsilon}{ t^+}-e+1\right)\right)+
\exp\left(- \epsilon \log\left(4\frac{\rho(y)+\epsilon}{
t^+}-e+1\right) \right)-2.
\end{align*}
Using (\ref{e:lotay}) we get
\begin{align*}
e^b+e^{-b}-2\le&  \epsilon^2 \left\{\left(4\frac{\rho(y)+\epsilon}{
t^+}-e+1\right)+\left(4\frac{\rho(y)+\epsilon}{
t^+}-e+1\right)^{-1}-2\right\}\\
\le& \epsilon^2 \left\{4\frac{\rho(y)+\epsilon}{ t^+}-e\right\},
\end{align*}
and hence
\begin{equation}\label{e:epsilon}
\frac{|h(t,z)-h(t,y)|^2}{4h(t,z)h(t,y)}\le \epsilon^2\left(
\frac{\rho(y)}{ t^+}+\frac{\epsilon}{t^+}-\frac{e}{4}\right)\le
\epsilon^2\left( \frac{\rho(y)}{
t^+}+\frac{1}{\tau}-\frac{e}{4}\right).
\end{equation}

Case III: $\rho(y)\wedge \rho(z)< 4^{-1}et^+<\rho(y)\vee \rho(z)$.
Since $|\rho(y)-\rho(z)|\le \epsilon$, we have $$
\rho(y)+\epsilon\ge \rho(y)\vee \rho(z)> 4^{-1}et^+~~{\rm
and~~}\rho(y)\vee\rho(z)-4^{-1}et^+< \epsilon.
$$ It implies
\begin{align*}
4\frac{\rho(y)+\epsilon}{t^+}-\frac{\rho(y)+\epsilon}{4^{-1}et^+}=\frac{\rho(y)+\epsilon}{4^{-1}et^+}\left(e-1\right)\ge
e-1.
\end{align*}
Hence
\begin{align*}
b=&
\left|(\rho(z)\vee \rho(y)-4^{-1}et^+)\log\left(\frac{\rho(z)\vee \rho(y)}{4^{-1}et^+}\right)\right|\\
\le&
\epsilon\log\left(\frac{\rho(y)+\epsilon}{4^{-1}et^+}\right)\le\epsilon\log\left(4\frac{\rho(y)+\epsilon}{t^+}-e+1\right).
\end{align*}
Similarly, we have (\ref{e:epsilon}) for this case.

On the other hand, note that $h(\cdot,y)$ is   differentiable on
$\bR^+$ and satisfies
\begin{align*}
-\frac{\partial }{\partial t}\log h(t,y)=&-\frac{\partial }{\partial
t}\left((\rho(y)-4^{-1}et^+)\log\left(1\vee\frac{\rho(y)}{4^{-1}et^+}\right)-\frac{t}{\tau}\right)\\
=&\frac{1}{\tau}+
4^{-1}e\log\left(1\vee\frac{\rho(y)}{4^{-1}et^+}\right)+\frac{
\big(\rho(y)-4^{-1}et^+\big)\vee 0}{t^+}\\
\ge& \frac{1}{\tau}+\left( \frac{\rho(y)}{
t^+}-\frac{e}{4}\right)\vee 0.
\end{align*}
 Therefore, in any case we have
 $$\frac{|h(t,z)-h(t,y)|^2}{4h(t,z)h(t,y)}\le \epsilon^2\left( \frac{1}{\tau}+\left( \frac{\rho(y)}{
t^+}-\frac{e}{4}\right)\vee 0\right)\le
-\epsilon^2\frac{\partial}{\partial t}\log h(t,y),$$
which implies $h\in \mathcal{F}(\bR_+)$.\qed\\
\\
The following two examples can be obtained in a similar way as Lemma
\ref{l:out} and we  leave it to the reader. See the examples in
\cite[Proposition 2.5 and Theorem 4.1]{CGZ} for a reference.
 \begin{example}\label{exa:1} Fix $a\in [0,\frac{1}{4}]$. Let
$h_1(t,x)=e^{a\rho(x)-\frac{a^2}{2}t}$. Then $h_1\in
\mathcal{F}(\bR_+)$.
\end{example}

\begin{example}\label{exa:2}  Fix $D\ge 5$, $R\ge 1$, $\Delta\ge\frac{24 R}{ D}$ and $s>
0$. For each $t\in [0,s]$ and $x\in \VV$, set
$h_2(t,x)=\exp\left(-\frac{\rho(x)^2}{D(s-t+\Delta)}\right)$. If
$1\le \rho(x)\le R$ for each $x\in\VV$,   then $h_2\in
\mathcal{F}([0,s])$.
\end{example}

Now, fix $o\in \VV$ and  for each  $R\ge 0$  set
\begin{align*}
\mathcal{G}_R(\bI)=& \{g: g \hbox{~is a function on~}
\bI\times\bR_+, ~g(t,r) \hbox{~is non-decreasing in
}r,~g(\cdot,d_{\nu}(o,\cdot)\wedge R) \in \mathcal{F}(\bI)\}.
\end{align*}
For brevity, we write $B_R=\{z\in \VV:d_\nu(o,z)<R\}$. The  lemma
below shows the way we use Corollary \ref{c:2.2*}.

\begin{lemma}\label{l:keyl}
Let $T\ge \tau\ge 0$ and $R\ge r\ge 0$. Let $u\in
\mathcal{H}([\tau,T])$ and   $g\in \mathcal{G}_R([\tau,T])$. Then
$$
\la u(T,\cdot)^2, 1-1_{B_R}\ra \le \frac{g(\tau, r)}{g(T,R)}\|
u(\tau,\cdot)\|^2 +\frac{g(\tau,R)}{g(T,R)}\la u(\tau,\cdot)^2,
1-1_{B_r}\ra.
$$
\end{lemma}
\pf Let $\rho(z)=\min\{d_\nu(o,z),R\}$ for each $z\in \VV$. Then
$\rho= R$ on $\VV\setminus B_R$ and hence
$$
\la u(T,\cdot)^2, 1-1_{B_R}\ra \le \la~ u(T,\cdot)^2,
g(T,\rho(\cdot))~\ra~ g(T,R)^{-1}.
$$
 By Corollary \ref{c:2.2*} and the hypothesis $u\in \mathcal{H}([\tau,T])$ and
$g(\cdot, \rho(\cdot))\in \mathcal{F}([\tau,T])$, we have
$$
 \la~ u(T,\cdot)^2,~ g(T,\rho)~\ra\le \la~ u(\tau,\cdot)^2,~ g(\tau,\rho )~\ra.
$$
Using the  condition  that $g(t,\cdot)$ is a non-decreasing
function, we get
\begin{align*}
\la u(\tau,\cdot)^2, g(\tau,\rho)\ra\le& \la u(\tau,\cdot)^2,
1_{B_r}\ra
g(\tau,r)+\la u(\tau,\cdot)^2, 1-1_{B_r}\ra g(\tau,R)\\
\le& g(\tau,r)\|u(\tau,\cdot)\|^2 +g(\tau,R)\la u(\tau,\cdot)^2,
1-1_{B_r}\ra,
\end{align*}
proving the lemma. \qed
\\

Furthermore, we set
\begin{align*}
\mathcal{H}_o=&\{ u\in \mathcal{H}(\bR_+):
u(0,z)=\nu_o^{-1/2}1_{\{o\}}(z) {\rm~for~each~} z\in \VV\}.
\end{align*}

\begin{prop}\label{t:shorttime}Let $u\in \mathcal{H}_o$.
For any  $t,R> 0$,  we have
\begin{align*}
\la u(t,\cdot)^2, 1-1_{B_R}\ra\le
\begin{cases}
\exp\left(-\frac{R^2}{8t}\right) & \hbox{if~~}t\ge R,\\
 \exp\left(-R\log
\left(\frac{1.01 R}{t}\right)+120\right) & \hbox{if~~}t\le R.
\end{cases}
\end{align*}
\end{prop}
\pf   Consider $t\ge R$ first. Take $a=\frac{R}{4t}$ then $a\in
(0,\frac{1}{4}]$. For each $s\ge 0$ and $r\ge 0$, set
$$
g_1(s,r)= e^{ar-\frac{a^2}{2}s}.
$$
By  Example \ref{exa:1},    $g_1\in \mathcal{G}_R(\bR_+)$. Use Lemma
\ref{l:keyl} and get for $r\in (0,R]$,
$$
\la u(t,\cdot)^2, 1-1_{B_R}\ra \le  \frac{g_1(0, r)}{g_1(t,R)}\|
u(0,\cdot)\|^2 +\frac{g_1(0,R)}{g_1(t,R)}\la u(0,\cdot)^2,
1-1_{B_r}\ra.
$$
From  $   u(0,z)=\nu_o^{-1/2}1_{\{o\}}(z)$, it follows immediately
that
$$
\la u(0,\cdot)^2, 1-1_{B_r}\ra=0~~~~{\rm and}~~~~\|u(0,\cdot)\|^2=1.
$$
So,
$$
\la u(t,\cdot)^2, 1-1_{B_R}\ra \le \lim_{r\rightarrow
0+}\frac{g_1(0, r)}{g_1(t,R)}= \frac{g_1(0, 0)}{g_1(t,R)}.
$$
Obviously, $g_1(0,0)=1$ and hence
$$
\la u^2(t,\cdot), 1-1_{B_R}\ra\le  e^{-aR+\frac{a^2}{2}t}.
$$
Substituting the value of $a$ into the above, we get the first
inequality of the proposition.\\

Next, suppose $t\le R$.  Choose  $\tau=(4c/e-1)t$, where
$b=(4c/e-1)^{-1}\approx 117.6$ and $c=e^{-e^{-1}}/1.01$.
 For each
$s\ge 0$ and $r\ge 0$, set
$$g_2(s,r)=\exp\left\{\big(r-4^{-1}e~(s+\tau)\big)~\log\left(1\vee\frac{r}{4^{-1}e~(s+\tau)}\right)
-\frac{s}{\tau}\right\}.$$ Obviously, $g_2(0,0)=1$.  By Lemma
\ref{l:out}, we have $g_2\in \mathcal{G}_R(\bR_+)$.  Since $x\log
(R/x)\le e^{-1}R$ for any $x>0$, we get
\begin{align}
\nonumber \log \big(g_2(t,R)\big)=&(R-ct)\log \left(\frac{R}{ct}\right)-b\\
\nonumber \nonumber
=&R\log\left(\frac{1.01R}{t}\right)+R\log\left(\frac{1}{1.01c}\right)-ct\log
\left(\frac{R}{ct}\right)-b\\
\nonumber\label{e:gggg} \ge&
R\log\left(\frac{1.01R}{t}\right)+R\log\left(\frac{1}{1.01c}\right)-e^{-1}R-120  \\
=&R\log\left(\frac{1.01R}{t}\right)-120.
\end{align}
From (\ref{e:gggg}) and $g_2\in \mathcal{G}_R(\bR_+)$, we prove the
second inequality of the proposition in the same way  as we did the
first. \qed
\begin{cor}\label{t:0011}For any  $z\in \VV$,
\begin{align*}
\P_o(X_t=z)\le\begin{cases} (\nu_z/\nu_o)^{1/2}
\exp\left\{-\frac{r^2}{16t}\right\}& \hbox{if~}t\ge r> 0;\\
(\nu_z/\nu_o)^{1/2} \exp\left(-\frac{r}{2}\log
\left(\frac{1.01r}{t}\right)+60\right)& \hbox{if~}r\ge t>
0,\\
\end{cases}
\end{align*}
where $r={d_\nu}(o,z)$.
\end{cor}
\pf  Recall the definition  $u_B$ in (\ref{e:exiting}). Denote by
$d(\cdot,\cdot)$ the graph distance of $\Gamma$. Set
$S_n=\{z:d(o,z)<n\}$. Then $S_n$ is a finite set since $\Gamma$ is a
locally finite graph and hence $u_{S_n}\in \mathcal{H}_o$. Clearly,
$u_{S_n}$ converges pointwise to $u$ as $n$ tends to infinity even
if the process $X$ explodes in finite time, where
$$
u(t,z)=\frac{\nu_o^{1/2}}{\nu_z}\P_o(X_t=z).
$$
Let $r=d_{\nu}(o,z)$, then we have $\la
u_{S_n}(t,\cdot)^2,1-1_{B_r}\ra\ge u_{S_n}(t,z)^2\nu_{z}$. So,
$$
u(t,z)^2\nu_z= \lim_{n\rightarrow \infty}u_{S_n}(t,z)^2\nu_z\le
\sup_n \la u_{S_n}(t,\cdot)^2,1-1_{B_r}\ra.
$$
Combining the above inequality with  Proposition \ref{t:shorttime},
we get the desired result. \qed

The long range bounds for transition probabilities were already
obtained by \cite[Theorems 2.1 and 2.2]{FL}, however, Corollary
\ref{t:0011} is more effective when $t\in [0.9r,1.1r]$    and
$r=d_{\nu}(o,z)$ is large. 

\section{Regular functions and integral estimates}
Recall that $A\ge 1$ and $\gamma>1$. Fix  $\delta\ge 1$,
$\theta_1=10^{-6}$ and $\theta_2=\theta_1/5$. Set
$$
\al=\min\{(2\gamma)^{-1}, (64\delta)^{-1}\}~~~~{\rm
and~~~~}\beta=\lceil \tfrac{\log 2}{\log \gamma}\rceil.
$$
 Let $u\in \mathcal{H}_o$ and $f:\bR_+\mapsto
\bR_+$ such that
\begin{equation}\label{e:subsolutiondiag}
\|u(t,\cdot)\|^2\le \frac{1}{f(2t)}~~~ \hbox{~for all~} t\in \bR_+.
\end{equation}
  In this section, we shall extend Proposition \ref{t:shorttime} into a result which can be used to prove Theorem \ref{t:main4}.

\begin{prop}\label{l:good}Suppose   $f$ is $(A,\gamma)-$regular and satisfies $f(t)\le A e^{\delta t}\hbox{~ for all~} t\in
\bR_+$. Then there exists  a constant $C_1>0$ which is independent
of $A,\gamma$ and $\delta$, such that  for  $t>0$,
\begin{equation}\label{e:good}
\Big\langle u(t,\cdot)^2, ~\exp\left(\theta_2\frac{
\big(d_{\nu}(o,\cdot)\wedge (2t)\big)^2 }{ t} \right)\Big\rangle\le
\frac{C_1 A^\beta}{f(2\alpha t)}.
\end{equation}
\end{prop}
 Before proving the proposition, we establish  some lemmas.

\begin{lemma}\label{l:highregular} If  $f$  is  an $(A,\gamma)-$regular function,  then
$$
f(2^{-k}t)\ge
\left(A^\beta\frac{f(t)}{f(\gamma^{-\beta}t)}\right)^{-k}
f(t)~~~~{\rm for~all~}k\in \bN~{\rm and~} t>0.
$$
\end{lemma}
\pf By the regularity, for any $t\ge s>0$ we have
\begin{align}
\label{e:regu512}\frac{f(\gamma^{\beta}
s)}{f(s)}=\prod_{j=0}^{{\beta}-1}
\frac{f(\gamma^{j+1}s)}{f(\gamma^js)}\le \prod_{j=0}^{{\beta}-1}
\left(A\frac{f(\gamma^{j+1}t)}{f(\gamma^{j}t)}\right)=A^{\beta}\frac{f(\gamma^{{\beta}}t)}{f(t)}.
\end{align}
In other words, an $(A,\gamma)$-regular function is also
$(A^\beta,\gamma^\beta)$-regular. Furthermore, by the monotonicity
we get
\begin{align*}
\frac{f(t)}{f(2^{-k}t)}\le \frac{f(t)}{f( \gamma^{-{\beta}k} t)}
=\prod_{j=-k}^{-1} \frac{f(\gamma^{\beta (j+1)}t )}{f(\gamma^{\beta
j}t)}\le \prod_{j=-k}^{-1}
\left(A^\beta\frac{f(t)}{f(\gamma^{-\beta}t)}\right)=\left(A^\beta\frac{f(t)}{f(\gamma^{-\beta}t)}\right)^k.
\end{align*}
\qed\\

\begin{lemma}\label{p:allneed}   If $f(t)\le A e^{\delta t}\hbox{~ for each~} t\in
\bR_+$, then there exists a constant  $c>0$ which is independent of
$A$ and $\delta$,  such that
\begin{align*}
\la u(t,\cdot)^2, 1-1_{B_R}\ra\le \frac{cA}{f(\frac{
R}{32\delta})}e^{-10^{-4}R}  {\rm~~~~for~}t>0{\rm~and~}R\in [t,64t].
\end{align*}
\end{lemma}
\pf Fix  $t>0$, $R\in [t,64 t]$, $x=t/R$ and $a=(64\delta)^{-1}$.
Then $a\le x \le 1$. Write $a_1=4^{-1}e~(a+0.45)$ and
$b=4^{-1}e(x+0.45)$. Then,  $$a_1\ge 4^{-1}e\cdot 0.45\ge 0.3 {\rm~~
and~~~}a_1\le  b\le 4^{-1}e(1+0.45)\le 0.99.$$  For each $s\ge 0$
and $r\ge 0$, we define
$$g(s,r)=\exp\left\{\big(r-4^{-1}e~(s+0.45R)\big)~\log\left(1\vee\frac{r}{4^{-1}e~(s+0.45R)}\right)
-\frac{s}{0.45R}\right\}.$$ By Lemma \ref{l:out},  we have $g\in
\mathcal{G}_R(\bR_+)$.  Applying  Lemma \ref{l:keyl} gives
\begin{equation}\label{e:xp}
\la u(t,\cdot)^2, 1-1_{B_R}\ra\le \frac{g(aR, a_1R)}{g(xR,R)}\|
u(aR,\cdot)\|^2 +\frac{g(aR,R)}{g(xR,R)}\la u(aR,\cdot)^2,
1-1_{B_{a_1R}}\ra.
\end{equation}
By a direct calculation, we get $g(aR,a_1R)\le 1$,
\begin{align*}
\log\left(g(aR,R)\right)\le &R(1-a_1)\log
\left(1\vee\frac{1}{a_1}\right)\le R(1-0.3)\log
\left(\frac{1}{0.3}\right)\le 0.8428R,
\end{align*}
and
\begin{align*}
\log\left(g(xR,R)\right)= & R(1-b)\log
\left(\frac{1}{b}\right)-\frac{x}{0.45}\ge R(1-0.99)\log
\left(\frac{1}{0.99}\right)-3 \ge 0.0001R -3.
\end{align*}
Thus, (\ref{e:xp}) becomes
\begin{align*}
\la u(t,\cdot)^2, 1-1_{B_R}\ra\le& e^{-0.0001R+3}\| u(aR,\cdot)\|^2
+e^{0.843R+3}\la u(aR,\cdot)^2, 1-1_{B_{a_1R}}\ra.
\end{align*}
By  (\ref{e:subsolutiondiag})  and the hypothesis $f(s)\le
Ae^{\delta s}$, we obtain,
\begin{align*}
\la u(t,\cdot)^2, 1-1_{B_R}\ra\le \frac{1}{f(2aR)}e^{-0.0001R+3}
+\frac{Ae^{2a\delta R}}{f(2aR)}e^{0.843R+3}\la u(aR,\cdot)^2,
1-1_{B_{a_1R}}\ra.
\end{align*}
By Proposition \ref{t:shorttime},
$$
\la u(aR,\cdot)^2, 1-1_{B_{a_1R}}\ra\le
\exp\left(-a_1R\log\left(\frac{a_1}{a} \right)+120\right).
$$
Therefore,
\begin{align*}
\la u(t,\cdot)^2, 1-1_{B_{R}}\ra\le& \frac{Ae^{123}}{f(2aR)}\left(
e^{-0.0001R}+\exp\left(2a\delta R+
0.843R-a_1R\log\left(\frac{a_1}{a} \right) \right)~\right).
\end{align*}
Substitute  $a=(64\delta)^{-1}$ and  get,
\begin{align*}
\la u(t,\cdot)^2, 1-1_{B_{R}}\ra\le&
\frac{Ae^{123}}{f(\frac{R}{32\delta})}\left(
e^{-0.0001R}+e^{-RC}~\right),
\end{align*}
 where $ C=a_1\log\left( 64a_1 \delta \right)-0.8743$.  Since
 $a_1\ge 0.3$ and $\delta\ge 1$, we have $
C\ge 0.01.$ So,
$$
\la u(t,\cdot)^2, 1-1_{B_R}\ra\le
\frac{2Ae^{123}}{f(\frac{R}{32\delta})}e^{-0.0001R}.
$$
\qed\\

\begin{prop}\label{p:favor} Suppose that  $f$ is $(A,\gamma)-$regular and satisfies $f(t)\le A e^{\delta t}\hbox{~ for all~} t\in
\bR_+$. Then there exists  a constant $C_0>0$ which is independent
of $A,\gamma$ and $\delta$, such that
$$
\la u(t,\cdot)^2, 1-1_{B_R}\ra\le \frac{ C_0A^\beta}{f(2\alpha
t)}\exp\left(-{\theta_1} \frac{R^2}{t}\right)~~{\rm~~for~all~}t\ge
R\ge 10^3.
$$
\end{prop}
\pf  Fix
$L=\log\left(A^\beta\frac{f(2t)}{f(2t/\gamma^\beta)}\right)$,
 $D=100$ and
$\Delta=\frac{R}{4}$.   If ${\theta_1}\frac{R^2}{t}-
L-\frac{1}{D\Delta}< {\theta_1}$, then we complete the proof since
\begin{align*}
\la u(t,\cdot)^2, 1-1_{B_R}\ra\le& \|u(t,\cdot)\|^2\\
\le&\frac{1}{f(2t)} \le\frac{e^{{\theta_1}}}{
f(2t)}\exp\left(L+\frac{1}{D\Delta}-{\theta_1}
\frac{R^2}{t}\right)\\
=&\frac{e^{{\theta_1}} A^\beta\exp({\frac{1}{D\Delta}})}{f(2
t/\gamma^\beta)}\exp\left(-{\theta_1}\frac{R^2}{t}\right)\\
\le& \frac{e^{{\theta_1}} A^\beta \exp(1/100)}{f(
t/\gamma)}\exp\left(-{\theta_1}\frac{R^2}{t}\right),
\end{align*}
where the last inequality uses the monotonicity of $f$. Therefore,
we may assume that
\begin{equation}\label{e:assume1}t\ge R\ge 10^3~~~~{\rm and~~~}{\theta_1}\frac{R^2}{t}-L-\frac{1}{D\Delta}\ge {\theta_1}.\end{equation}
This implies that $R\le t\le  R^2$ and $L\le {\theta_1}\frac{R^2}{t}$.\\

Let $\rho(x)=(R-d_\nu(o,x))\vee 1$ for any $x\in \VV$. Then $\rho$
satisfies (\ref{e:rholimitation1}) and $1\le \rho(x)\le R$. For each
$s\in [0,t]$ and $r\ge 0$, set
$$
g(s,r)=\exp\left(-\frac{\big((R-r)\vee
1\big)^2}{D(t-s+\Delta)}\right).
$$
Then $g\in \mathcal{G}_R([0,t])$ by Example \ref{exa:2} and the
 argument above about $\rho$. From  Lemma \ref{l:keyl}, we get that
for any $r\in [0,R]$ and $s\in [0,t]$,
\begin{align}
\nonumber\la u(t,\cdot)^2& , 1-1_{B_R}\ra\le\frac{g(s, r)}{g(t,R)}\|
u(s,\cdot)\|^2 +\frac{g(s,R)}{g(t,R)}\la u(s,\cdot)^2,
1-1_{B_{r}}\ra~~~~~~~\\
\label{e:vip}\le& \frac{\exp({\frac{1}{D\Delta}})}{f(2s)}
\exp\left(-\frac{(R-r)^2}{D(t-s+\Delta)}\right)+\exp({\frac{1}{D\Delta}})\la
u(s,\cdot)^2, 1-1_{B_{r}}\ra.
\end{align}
We shall iterate using (\ref{e:vip}).
 Let us build a sequence $\{(t_j, R_j): 0\le j\le j_0\}$.  Take
$$
 t_j=t/2^{j-1},~~R_j=R/2+ R/(j+1)  \hbox{~for each~}0\le j\le j_0;
$$
and  $$j_0=\min\{j: R_j\ge t_j\}.$$ Then $j_0\ge 1$ and for all
$0\le j<j_0$ we have $t_j>R_j>R/2>1$. Hence
$$
t_j-t_{j+1}=t_j/2\ge R/4=\Delta.
$$
From $t_{j_0-1}>R/2$, we  get
$$
j_0< \frac{\log (8t/R)}{\log 2}.
$$
Using the identity $
\left(R_j-R_{j+1}\right)^2=\frac{R^2}{(j+1)^2(j+2)^2},$ we obtain
$$
\frac{(R_j-R_{j+1})^2}{D(t_j-t_{j+1}+\Delta)}\ge
\frac{(R_j-R_{j+1})^2}{D t_j}
  =\frac{2^{j-1}}{D(j+1)^2(j+2)^2 } \frac{R^2}{t}.
$$
Note that
$$
\min\left\{ \frac{2^{j-1}}{100(j+1)^3(j+2)^2 }: j\ge
1\right\}=\frac{2^{6-1}}{100(6+1)^3(6+2)^2}\approx 1.5\times
10^{-5}.
$$
Sicne   $\theta_1=10^{-6}$, it follows immediately that
$$
\frac{(R_j-R_{j+1})^2}{D(t_j-t_{j+1}+\Delta)}\ge    (j+1)
{\theta_1}\frac{R^2}{t}.
$$
 Iterating  (\ref{e:vip}), we obtain
\begin{align*}
&\la u(t,\cdot)^2, 1-1_{B_R}\ra=\la u(t_1,\cdot)^2, 1-1_{B_{R_1}}\ra\\
\le& \sum_{j=1}^{j_0-1}\frac{\exp(\frac{j}{D\Delta})}{
f(2t_{j+1})}\exp\left(-\frac{(R_j-R_{j+1})^2}{D(t_j-t_{j+1}+\Delta)}\right)+\exp\left(
\frac{j_0-1}{D\Delta}\right)\la u(t_{j_0},\cdot)^2, 1-1_{B_{R_{j_0}}}\ra\\
:=&\Lambda_1+\Lambda_2.
\end{align*}
 By Lemma \ref{l:highregular}, we have \begin{equation}\label{e:usedrr} f(2t_{j+1})\ge
f(2t)e^{-jL}.\end{equation} Using (\ref{e:assume1}), we conclude
\begin{align}
\nonumber\Lambda_1\le& \frac{1}{ f(2t)}\exp\left(-{\theta_1}
\frac{R^2}{t}\right)\sum_{j=1}^{j_0-1}\exp\left(-j\left({\theta_1}\frac{R^2}{t}-L-\frac{1}{D\Delta}\right)\right)\\
\nonumber\le& \frac{1}{ f(2t)}\exp\left(-{\theta_1}
\frac{R^2}{t}\right)\sum_{j=1}^{j_0-1}\exp(-j{\theta_1})\\
\label{e:Lambda1}\le&
\frac{e^{-{\theta_1}}(1-e^{-{\theta_1}})^{-1}}{
f(2t)}\exp\left(-{\theta_1} \frac{R^2}{t}\right).
\end{align}

On the other hand, since $2t_{j_0}=t_{j_0-1}>R_{j_0-1}>R_{j_0}\ge
t_{j_0}$, we use Lemma \ref{p:allneed} to get
\begin{equation}\label{e:usedexp}
\la u(t_{j_0},\cdot)^2, 1-1_{B_{R_{j_0}}}\ra\le
\frac{cA}{f(\frac{R_{j_0}}{32\delta})}e^{-10^{-4}R_{j_0}},
\end{equation}
where $c$ is a  constant which is independent of $A,\gamma$ and
$\delta$. By Lemma \ref{l:highregular} and (\ref{e:regu512}), we
also have
\begin{align}
\nonumber f(\frac{R_{j_0}}{32\delta})\ge& f(\frac{t_{j_0}}{32\delta})\ge \left(A^\beta\frac{f(\frac{t}{32\delta})}{f(\frac{t}{32\delta\gamma^\beta})}\right)^{-j_0+1}f(\frac{t}{32\delta})\\
\nonumber \ge& \left(A^{2\beta}\frac{f(2t)}{f(2t/\gamma^\beta)}\right)^{-j_0+1}f(\frac{t}{32\delta})\\
\label{e:usedreg} \ge& f(\frac{t}{32\delta})e^{-2j_0L}.
\end{align}
    So,
\begin{align*}
\Lambda_2=&\exp( \frac{j_0-1}{D\Delta})\la u(t_{j_0},\cdot)^2,
1-1_{B_{R_{j_0}}}\ra
\le\exp(\frac{j_0-1}{D\Delta})\frac{cA}{f(\frac{R_{j_0}}{32\delta})}e^{-10^{-4}R_{j_0}}\\
\le& \frac{cA}{f(\frac{t}{32\delta})}
\exp\left(\frac{j_0}{D\Delta}+2j_0L-10^{-4}R/2\right).
\end{align*}
Note that \\
$$10^3\le R\le t\le R^2,~~~~j_0< \frac{\log( 8t/ R)}{\log 2},~~~~D\Delta=25R~~~~{\rm and~~~~}L\le {\theta_1}
\frac{R^2}{t}.$$ From these inequalities, we calculate
\begin{align*}
\frac{j_0}{D\Delta R}\le& \frac{\log(8t/R)}{25R^2\log 2}\le
\frac{\log(8R)}{25R^2\log 2}\le \frac{\log(8\cdot 10^3)}{25\cdot
10^6\cdot\log 2}< 5.2\times 10^{-7};\\
\frac{2j_0L}{R}\le& 2{\theta_1}\frac{\log( 8t/ R)}{\log 2}
\frac{R}{t}\le 2\theta_1\frac{8}{e\log 2}< 8.5\times 10^{-6}.
\end{align*}
So,  $\frac{j_0}{D\Delta}+2j_0L-10^{-4}R/2< -\theta_1R$ and hence
\begin{equation}\label{e:22.1}
\Lambda_2\le \frac{cA}{ f(\frac{t}{32\delta})}e^{-\theta_1R}\le
\frac{cA}{ f(\frac{t}{32\delta})}e^{-\theta_1R^2/t}.
\end{equation}

Finally, we choose
$$
C_0=e^{\theta_1+0.01} +e^{-{\theta_1}}(1-e^{-{\theta_1}})^{-1}+c
$$
and
  complete the proof.\qed\\

\noindent {\it \bf  Proof of Proposition \ref{l:good}.}  Write
$\rho(z)=d_{\nu}(o,z)\wedge (2t)$ for short. If $t\le 10^6$, then
the result is trivial since
$$
\Big\langle u(t,\cdot)^2, ~\exp\left({\theta_2}\frac{ \rho^2 }{ t}
\right)\Big\rangle\le e^{4{\theta_2} t}\|u(t,\cdot)\|^2\le
\frac{e^{4\cdot10^6{\theta_2}}}{f(2t)}.
$$
So, we may assume that $t\ge 10^6$ in the following.

Fix $R=t^{1/2}$ and $n=\lceil\frac{\log (t/R)}{\log 2}\rceil$. Then
$2^nR\ge t$, and $t\ge 2^{j-1}R\ge 10^3$ for each $1\le j\le n$.
Write
$$\Upsilon_0=\la u(t,\cdot)^2,~ e^{{\theta_2} \rho^2/ t}1_{B_R}\ra,
~\Upsilon_\infty=\la u(t,\cdot)^2,~~ e^{{\theta_2} \rho^2/
t}(1-1_{B_{t}})\ra  $$ and set
$$
\Upsilon_j=\la u(t,\cdot)^2,~ e^{{\theta_2} \rho^2/
t}(1_{B_{2^{j}R}}-1_{B_{2^{j-1}R}})\ra~~~{\rm for~}1\le j\le n.
$$
Then
$$
\Big\langle u(t,\cdot)^2, ~\exp\left({\theta_2}\frac{ \rho^2 }{ t}
\right)\Big\rangle\le
\Upsilon_0+\sum_{j=1}^n\Upsilon_j+\Upsilon_\infty.
$$
We estimate each $\Upsilon_j$ separately.

The first term  admits the estimate
$$
\Upsilon_0\le \la u(t,\cdot)^2, e^{{\theta_2} }1_{B_R}\ra \le
e^{\theta_2}\|u(t,\cdot)\|^2\le \frac{e^{\theta_2}}{ f(2t)}.
$$
Next, for each $1\le j\le n$,   we have
\begin{equation}\label{e:larupsi}
\Upsilon_j\le \la u(t,\cdot)^2,~ e^{{\theta_2}
(2^j)^2}(1_{B_{2^{j}R}}-1_{B_{2^{j-1}R}})\ra\le  e^{
4^j{\theta_2}}~\la u(t,\cdot)^2,~ 1-1_{B_{2^{j-1}R}}\ra .
\end{equation}
 Set $C_0$  as in
Proposition \ref{p:favor}. Then
$$
\la u(t,\cdot)^2,~ 1-1_{B_{2^{j-1}R}}\ra\le \frac{C_0A^\beta}{
f(2\alpha t)}\exp\left(-{\theta_1}\cdot4^{j-1}\right).
$$
 By definition ${\theta_2}={\theta_1}/5$; therefore we get
\begin{align*}
\Upsilon_j\le\frac{C_0 A^\beta }{ f(2\alpha
t)}\exp\left(-{\theta_2}\cdot4^{j-1}\right).
\end{align*}
For the remaining term,
$$
\Upsilon_\infty\le  e^{4{\theta_2} t}  \la u(t,\cdot)^2,~~
(1-1_{B_{t}})\ra.
$$
Using Proposition \ref{p:favor} again gives
$$
\Upsilon_\infty\le e^{4{\theta_2} t}\cdot
\frac{C_0A^\beta}{f(2\alpha t)}e^{-{\theta_1}
t}=\frac{C_0A^\beta}{f(2\alpha t)}e^{-{\theta_2} t}\le
\frac{C_0A^\beta}{f(2\alpha t)}.
$$

Therefore,
\begin{align*}
\Big\langle u(t,\cdot)^2, ~\exp\left({\theta_2}\frac{ \rho^2 }{ t}
\right)\Big\rangle\le& \frac{e^{\theta_2}}{ f(2t)}+\sum_{j=1}^n
\frac{C_0 A^\beta }{ f(2\alpha
t)}\exp\left(-{\theta_2}\cdot4^{j-1}\right)
+\frac{C_0A^\beta}{f(2\alpha t)}\\
\le&\frac{C_1A^\beta}{f(2\alpha t)},
\end{align*}
where
$$
C_1=e^{4\cdot10^6{\theta_2}}+
C_0\sum_{j=1}^\infty\exp\left(-{\theta_2}\cdot4^{j-1}\right)+C_0.
$$
\qed
\section{Proof of Theorem \ref{t:main4}}
 {\it\bf Proof of Theorem \ref{t:main4}.} Recall the notation
 $\mathcal{H}_{\cdot}$ from Section 2.
 Fix $ t\ge  {d_\nu}(x_1,x_2)$ and
$s=t/2$.  For each $z\in \VV$ and $i\in\{1,2\}$, set
$$\rho_i(s,z)=d_\nu(x_i,z)\wedge (2 s){\rm~~~~and~~~~}h_i(s,z)=\exp\left(
\frac{1}{2}\cdot{\theta_2}\frac{\rho_i(s,z)^2}{ s}\right).$$ Then
$2\rho_1(s,z)^2+2\rho_2(s,z)^2\ge d_\nu(x_1,x_2)^2$ and so
\begin{equation}\label{e:hequa}
h_1(s,z)h_2(s,z)\ge\exp\left(\frac{{\theta_2}}{2}\cdot\frac{d_{\nu}(x_1,x_2)^2}{
t}\right) .
\end{equation}
Let $d(\cdot,\cdot)$ be the graph distance of $\Gamma$. As in
Corollary \ref{t:0011}, we define
$$
u_{ij}(s,z)=\frac{\nu_{x_i}^{1/2}}{\nu_z}\P_{x_i}(X_s=z,~~
\inf\{l\in \bR_+: d(x_i,X_l)\ge j\}>s)
$$
and  $u_i(s,z)=\frac{\nu_{x_i}^{1/2}}{\nu_z}\P_{x_i}(X_s=z)$. Then
$\{u_{ij}(s,z): j=1,2,\cdots\}$ is a non-decreasing sequence and
satisfies
$$
\lim_{j\rightarrow \infty} u_{ij}(s,z)=u_i(s,z).
$$
By (\ref{e:diag}), for any $l\ge 0$ we have
$$
\| u_{ij}(l,\cdot)\|^2\le \|
u_{i}(l,\cdot)\|^2=\P_{x_i}(X_{2l}=x_i)\le \frac{1}{f_i(2l)}.
$$
Since  $u_{ij}\in \mathcal{H}_{x_i}$, we use  Proposition
\ref{l:good} and get
$$
\|u_{ij}(s,\cdot) h_i(s,\cdot)\|^2=\left\langle u_{ij}(s,\cdot)^2,~
\exp\left({\theta_2}\frac{\rho_i(s,\cdot)^2}{ s}\right)\right\rangle
\le \frac{C_1A^\beta}{f_i(2\alpha s)}=\frac{C_1A^\beta}{f_i(\alpha
t)} .
$$
By the Monotone Convergence Theorem,
$$
\|u_i(s,\cdot)
h_i(s,\cdot)\|^2=\lim_{j\rightarrow\infty}\|u_{ij}(s,\cdot)
h_i(s,\cdot)\|^2  \le \frac{C_1A^\beta}{f_i(\alpha t)} .
$$
 By (\ref{e:hequa}) and the Cauchy-Schwarz inequality,
we obtain
\begin{align*}
\P_{x_1}(X_{t}=x_2)=&\sum_{z\in\VV}\P_{x_1}(X_{s}=z)\P_{x_2}(X_{s}=z)\frac{\nu_{x_2}}{\nu_z}\\
=& (\nu_{x_2}/\nu_{x_1})^{1/2} ~\la u_1(s,\cdot), u_2(s,\cdot)\ra\\
\le &(\nu_{x_2}/\nu_{x_1})^{1/2}~ \big\langle~ u_1(s,\cdot)
h_1(s,\cdot),~ u_2(s,\cdot) h_2(s,\cdot)~\big\rangle ~
\exp\left(-\frac{{\theta_2}}{2}\cdot\frac{d_{\nu}(x_1,x_2)^2}{
t}\right) \\
 \le&(\nu_{x_2}/\nu_{x_1})^{1/2}~ \| u_1(s,\cdot) h_1(s,\cdot)\|~~\|
u_2(s,\cdot) h_2(s,\cdot)\| ~
\exp\left(-\frac{{\theta_2}}{2}\cdot\frac{d_{\nu}(x_1,x_2)^2}{
t}\right) \\
\le&\frac{C_1A^\beta (\nu_{x_2}/\nu_{x_1})^{1/2}}{\sqrt{f_1(\alpha
t)f_2(\alpha t)}}~
\exp\left(-\frac{{\theta_2}}{2}\cdot\frac{d_{\nu}(x_1,x_2)^2}{
t}\right).
\end{align*}
Set $\theta={\theta_2}/2$ and we complete the proof. \qed\\

\section{Regularity on an interval}
{\it\bf Proof of Theorem \ref{t:main5}.} First, we show that
\begin{equation}\label{e:00940}
\Big\langle u(t,\cdot)^2, ~\exp\left(\theta_2\frac{
\big(d_{\nu}(o,\cdot)\wedge (2t)\big)^2 }{ t} \right)\Big\rangle\le
\frac{C_1 A^\beta}{f(2\alpha t)}~~~{\rm for~~~}t\in
[(2\alpha^{-1}T_1)^2,T_2/2).
\end{equation}
Take $t\in [(2\alpha^{-1}T_1)^2,T_2/2)$, $t_j=t/2^{j-1}$,
$R_j=R/2+R/(j-1)$, $j_0=\min\{j:R_j\ge t_j\}$ and $R=t^{1/2}$ as in
Propositions \ref{l:good} and \ref{p:favor}. Then $$ 2\al t_{j_0}>
\al R/2=\al t^{1/2}/2\ge T_1.
$$
Using (\ref{e:11181}) and the regular condition  on $[T_1,T_2)$, we
 still have the inequalities (\ref{e:usedrr}), (\ref{e:usedexp}) and
(\ref{e:usedreg}). Therefore, we can get (\ref{e:00940}) in the same
way as we did Proposition \ref{l:good}. Furthermore, by
(\ref{e:00940}) and the Cauchy-Schwarz inequality,
 we finish the proof of Theorem \ref{t:main5} similar as that of  Theorem \ref{t:main4}. \qed\\

The lower bound of $\widetilde{T}_1$ in Theorem \ref{t:main5} can be
improved if one knows more information about the growth rate of
$f_i$. Heat kernels having either sub-exponential decay or
polynomial decay appear in many groups, see Hebisch and Saloff-Coste
\cite[Theorem 4.1]{HS}. More importantly, there are a lot of papers
which studied random walks on $\mathbb{Z}^d$ with random
conductances and showed that $p_t(x,x)\le c t^{-d/2}$ under
different conditions, such as [1-2] and [4-6]. Therefore, we have to
consider function
$f_i$ whose growth rate has either sub-exponential or polynomial.\\

Fix $A\ge 1,\gamma>1, \theta_1=10^{-6}$ and $
\theta=\theta_2/2=\theta_1/10$  as before.

\begin{thm}\label{t:main2} Let $\delta\ge 0$ and  $\epsilon\in [0,1)$. If each $f_i$ is  $(A,\gamma)$-regular
 on $[T_1,T_2)$ and satisfies
\begin{equation}\label{e:speed222} f_i(t)\le
Ae^{\delta t^{\epsilon}} {\rm ~for~ all~}t\in [T_1,T_2),
\end{equation}
then there exists a constant $C_1(A,\gamma,\delta,\epsilon)>0$ such
that for each $t\in[\widetilde{T_1} ,T_2)$,
\begin{equation}\label{e:intervalresult}
\P_{x_1}(X_{t}=x_2)\le \frac{C_1 (\nu_{x_2}/\nu_{x_1})^{1/2}
}{\sqrt{f_1(\frac{t}{2\gamma})f_2(\frac{t}{2\gamma})}}\exp\left(-\theta\frac{d_{\nu}(x_1,x_2)^2}{t}\right),
\end{equation}
where $\widetilde{T_1}=(2^9\delta~T_1^{1+\epsilon})\vee
d_{\nu}(x_1,x_2)$.
\end{thm}
\begin{thm}\label{t:main3} Let $\epsilon\ge 0$. If each $f_i$ is $(A,\gamma)-$regular on $[T_1,T_2)$ and
satisfies $$ f_i(t)\le At^{\epsilon}{\rm ~for~ all~}t\in [T_1,T_2),
$$
then there exists  a constant  $C_1(A,\gamma,\epsilon)>0$  such that
(\ref{e:intervalresult}) holds  for $t\in[\widetilde{T_1},T_2)$.
Here,
$$\widetilde{T}_1=\big(~2^{10}\epsilon~T_1\log (T_1\vee1)~\big)\vee
d_\nu(x_1,x_2).$$
\\
\end{thm}

Let's begin with  Theorem \ref{t:main2}.  As the proof of Theorem
\ref{t:main4}, we need some results which are similar to
Propositions \ref{p:favor} and \ref{l:good}.
\begin{prop}\label{p:favor3} Let $\delta> 0$ and  $\epsilon\in (0,1)$.  Let $u,f$ be defined as in Section 3.  Suppose
further that   $f$ is  $(A,\gamma)$-regular  on $[T_1,T_2)$ and
satisfies
\begin{equation}\label{e:speed222} f(t)\le
Ae^{\delta t^{\epsilon}} {\rm ~for~ all~}t\in [T_1,T_2).
\end{equation} Then there exists a constant
$C_0(A,\gamma,\delta,\epsilon)>0$
 such that  for $R\ge
\max\{4, 2\kappa^{\frac{1+\epsilon}{1-\epsilon}},~ 2(\kappa
T_1)^{(1+\epsilon)/2}\}$ and $t\in [\kappa^{-1}R^{2/(1+\epsilon)},
T_2/2)$, we have
$$
\la u(t,\cdot)^2, 1-1_{B_R}\ra\le \frac{
C_0}{f(\frac{t}{\gamma})}\exp\left(-{\theta_1} \frac{R^2}{t}\right).
$$
Here, $\kappa=(64\delta)^{1/(1+\epsilon)}$.
\end{prop}

\pf We only show the part of the proof which is  different from that
of Proposition \ref{p:favor}.

Fix  $R\ge \max\{4, 2\kappa^{\frac{1+\epsilon}{1-\epsilon}},~
2(\kappa T_1)^{(1+\epsilon)/2}\}$ and $t\in
[\kappa^{-1}R^{2/(1+\epsilon)}, T_2/2)$.
   Take
$L,D,\Delta, R_j$ and $t_j$ as in Proposition \ref{p:favor}. We may
still assume that $\theta_1\frac{R^2}{t}-L-\frac{1}{D\Delta}\ge
\theta_1$. (Hence $t\le R^2$ and $L\le \theta_1\frac{R^2}{t}$.)
However, we set
$$j_0=\min\{j: R_j^{2/(1+\epsilon)}\ge \kappa  t_j\}.$$
Since $R\ge \max\{4, 2\kappa^{(1+\epsilon)/(1-\epsilon)}\}$,  for
$j<j_0$ we have
$$t_j>\kappa^{-1}R_j^{2/(1+\epsilon)}>\kappa^{-1}(R/2)^{2/(1+\epsilon)}\ge R/2\ge 2.$$
 Hence
$$t_j-t_{j+1}=t_j/2\ge R/4=\Delta.$$  From $R\ge 2(\kappa
T_1)^{(1+\epsilon)/2}$, it deduces
$$t_{j_0}=t_{j_0-1}/2\ge \kappa^{-1}(R/2)^{2/(1+\epsilon)}/2\ge T_1/2.  $$
So,   $T_1\le 2t_{j+1}\le 2t< T_2$ for each $j<j_0$. By the
hypothesis that
   $f$ is $(A,\gamma)$-regular on $[T_1,T_2)$, one has
 $$f(2t_{j+1})\ge f(2t)e^{-jL} $$
 the same  as Lemma \ref{l:highregular}. Hence
(\ref{e:Lambda1}) holds under this circumstance, too. That is,
$$
\Lambda_1:=\sum_{j=1}^{j_0-1}\frac{\exp({\frac{j}{D\Delta}})}{
f(2t_{j+1})}\exp\left(-\frac{(R_j-R_{j+1})^2}{D(t_j-t_{j+1}+\Delta)}\right)\le
\frac{e^{-{\theta_1}}(1-e^{-{\theta_1}})^{-1}}{
f(2t)}\exp\left(-{\theta_1} \frac{R^2}{t}\right).
$$
\\

Next, if $t_{j_0}<R_{j_0}$ then by Proposition \ref{t:shorttime},
\begin{align*}
\la u(t_{j_0},\cdot)^2, 1-1_{B_{R_{j_0}}}\ra \le  c_1e^{-
2c_1^{-1}R_{j_0}}\le c_1e^{-c_1^{-1}R}\le c_2
\exp\left(-\frac{\kappa}{16} R^{2\epsilon/(1+\epsilon)}\right),
\end{align*}
where $c_1,c_2\ge 1$  are constants.
  If $t_{j_0}\ge R_{j_0}$ then by Proposition \ref{t:shorttime} we still have
\begin{align*}
\la u(t_{j_0},\cdot)^2, 1-1_{B_{R_{j_0}}}\ra \le
e^{-\frac{R_{j_0}^2}{8t_{j_0}}} \le \exp\left(-\frac{\kappa }{8}
R_{j_0}^{2\epsilon/(1+\epsilon)}\right) \le
c_2\exp\left(-\frac{\kappa}{16} R^{2\epsilon/(1+\epsilon)}\right).
\end{align*}
From $R_{j_0}^{2/(1+\epsilon)}\ge \kappa t_{j_0}$ and
$R_{j_0-1}^{2/(1+\epsilon)}< \kappa t_{j_0-1}$,  we get the
following  inequalities respectively:
$$t_{j_0}\le R^{2/(1+\epsilon)}/\kappa  {\rm~~~and~~~} j_0< \frac{1}{\log 2}\log \left( \frac{16\kappa t}{ R^{2/(1+\epsilon)}}  \right). $$
Hence
$$
f(2t)\le  f(2t_{j_0})e^{j_0L}\le
f(2R^{2/(1+\epsilon)}/\kappa)\exp\left\{ \frac{1}{\log 2} \log
\left( \frac{16\kappa t}{ R^{2/(1+\epsilon)}}  \right) \cdot
L\right\}
$$
By  (\ref{e:speed222}) and the assumption $L\le
\theta_1\frac{R^2}{t}$,
\begin{align*}
f(2t)\le& A\exp\left(\frac{ 2^\epsilon \delta}{\kappa^\epsilon}
R^{2\epsilon/(1+\epsilon)}\right)\cdot\exp\left\{ \frac{1}{\log 2}
\log \left( \frac{16\kappa t}{ R^{2/(1+\epsilon)}}  \right)\cdot  \theta_1\frac{R^2}{t}\right\} \\
\le&A\exp\left(\frac{2^\epsilon \delta}{\kappa^\epsilon}
R^{2\epsilon/(1+\epsilon)}\right)\exp\left\{ \frac{\theta_1}{e\log
2}
 16\kappa R^{2\epsilon/(1+\epsilon)} \right\}.
\end{align*}
Since $t\le R^2$ and $R\ge 4$, there exists a constant $c_3$  such
that
$$
\exp({\frac{j_0}{D\Delta}})\le  \exp\left( \frac{1}{\log 2}\log
\left( \frac{16\kappa t}{ R^{2/(1+\epsilon)}}
\right)\cdot\frac{1}{25R} \right)\le  \exp\left( \frac{1}{\log
2}\log \left( \frac{16\kappa R^2}{ R^{2/(1+\epsilon)}}
\right)\cdot\frac{1}{25R}  \right)\le c_3.
$$
Combining the above inequalities  gives
\begin{align*}
\Lambda_2:=&\exp({\frac{j_0-1}{D\Delta}})\la u(t_{j_0},\cdot)^2,
1-1_{B_{R_{j_0}}}\ra \\
\le& c_3c_2\exp\left(-\frac{\kappa}{16}
R^{2\epsilon/(1+\epsilon)}\right)\\
\le& c_3c_2 \exp\left(-\frac{\kappa}{16}
R^{2\epsilon/(1+\epsilon)}\right)\cdot \frac{1}{f(2t)}\cdot
A\exp\left(\frac{2^{\epsilon}\delta}{\kappa^\epsilon}
R^{2\epsilon/(1+\epsilon)}\right)\exp\left\{ \frac{\theta_1}{e\log
2}
 16\kappa R^{2\epsilon/(1+\epsilon)} \right\}\\
 =&\frac{c_4}{f(2t)}\exp\left(\left(-\frac{\kappa}{16}+\frac{2^\epsilon\delta}{\kappa^\epsilon}+\frac{\theta_1}{e\log
2}
 16\kappa
 \right)
R^{2\epsilon/(1+\epsilon)}\right),
\end{align*}
where $c_4=c_3c_2A$.
 Substituting  $\kappa=(64\delta)^{1/(1+\epsilon)}$ and using
the condition $ t\ge \kappa^{-1} R^{2/(1+\epsilon)}$,
\begin{align*}
\Lambda_2\le
&\frac{c_4}{f(2t)}\exp\left(\left(-\frac{\kappa}{16}+\frac{2^\epsilon\kappa}{64}+\frac{\theta_1}{e\log
2}
 16\kappa
 \right)
R^{2\epsilon/(1+\epsilon)}\right)\\
\le& \frac{c_4}{f(2t)}\exp\left(-\frac{\kappa}{64}
R^{2\epsilon/(1+\epsilon)}\right)\\
\le&\frac{c_4}{f(2t)}\exp\left(-\frac{1}{64} \frac{R^2}{t}\right).
\end{align*}
This completes  the proof.\qed

\begin{prop}\label{l:good3} Under the condition of Proposition \ref{p:favor3},  there
exists a constant   $C_0(A,\gamma,\delta,\epsilon)>0 $ such that
$$
\left\langle u(t,\cdot)^2,
\exp\left(\theta_2\frac{(d_\nu(o,\cdot)\wedge(2t))^2}{t}\right)\right\rangle\le
\frac{ C_0}{f(\frac{t}{\gamma})}~~~~{\rm~for~}t\in [2^8\delta
T_1^{1+\epsilon}, T_2/2).
$$
\end{prop}
\pf We only show the difference from Proposition \ref{l:good}. Fix
$\kappa=(64\delta)^{1/(1+\epsilon)}$ and $t_0=\max\{16,
4\kappa^{(2+2\epsilon)/(1-\epsilon)},
\kappa^{-(1+\epsilon)/\epsilon}\}$.  Let $t\in [2^8\delta
T_1^{1+\epsilon}, T_2/2)$.     If $t\le t_0$, then as before the
result is trivial. So, we may assume further $t\ge t_0$. Fix
$R=t^{1/2}$. Then $$R\ge \max\{4,
2\kappa^{\frac{1+\epsilon}{1-\epsilon}},~2(\kappa
T_1)^{(1+\epsilon)/2}\}   {\rm~~~~~and~~~~~}\kappa t\ge
R^{2/(1+\epsilon)}.$$ Define $\theta_2, \Upsilon_j$ and $n$ as in
Proposition \ref{l:good}. However, we set $$m=\max\{j: \kappa t\ge
(2^j R)^{2/(1+\epsilon)}\}.$$ Then by  Proposition \ref{p:favor3},
for $1\le j\le m\wedge n$ we have
$$
 \Upsilon_j\le \frac{C_0}{f(t/\gamma)}\exp\left(-\theta_2\cdot
 4^{j-1}\right).
$$
If $m+1\le j\le n$ then use Proposition \ref{t:shorttime} and get
\begin{align*}
\Upsilon_j\le&  e^{4^j \theta_2} \la u(t,\cdot)^2,
1-B_{2^{j-1}R}\ra\le e^{4^j \theta_2} \cdot
\exp\left(-\frac{(2^{j-1})^2}{8} \right)\le \exp\left(-\frac{
4^{j-1}}{12} \right).
\end{align*}
By the definition of $m$, one has
 $\kappa t<(2^{m+1}R)^{2/(1+\epsilon)}=(2^{m+1}t^{1/2})^{2/(1+\epsilon)} $ and so,
  \begin{equation}\label{e:meps}4^{m+1}>
 \kappa^{1+\epsilon}t^\epsilon=64\delta
 t^\epsilon.\end{equation}
By (\ref{e:speed222}) and (\ref{e:meps}), we still have
\begin{align*}
\Upsilon_j\le& \frac{Ae^{\delta t^\epsilon}}{f(t)}\exp\left(-\frac{
4^{j-1}}{12} \right)\le \frac{A}{f(t)}\exp\left(4^{m-2}-\frac{
4^{j-1}}{12}\right)\\
\le& \frac{A}{f(t)}\exp\left(4^{j-3}-\frac{
4^{j-1}}{12}\right)=\frac{A}{f(t)}\exp\left(-\frac{
4^{j-1}}{48}\right).
\end{align*}
For the other terms $\Upsilon_0$ and $\Upsilon_\infty$, one can get
the estimates  the same as we did in Proposition \ref{l:good} and so
we finish the
proof.\qed\\

\noindent {\it\bf Proof of Theorem \ref{t:main2}.} If
$\epsilon\delta=0$, then the problem is reduced to Corollary
\ref{t:0011} since each $f_i$ has a constant upper bound on
$[T_1,T_2)$. Otherwise, if $\delta>0$ and $\epsilon\in (0,1)$ then
we can get the proof  as Theorem \ref{t:main4} by  using Proposition
\ref{l:good3} and the Cauchy-Schwarz inequality.\qed

\noindent {\it\bf Proof of Theorem \ref{t:main3}.}  We obtain  a
similar result as Proposition \ref{p:favor3} just by setting
$$j_0=\min\{j: R_j^2/\log R_j\ge \kappa t_j\},$$ and then
prove the theorem  as above.  \qed\\

Enlightened  by  Boukhadra, Kumagai and  Mathieu \cite{BKM},  we
give an application of Theorem 5.1. Set
$$p_t(x,y)=\frac{\P_x(X_t=y)}{\nu_y}$$ for the heat kernel of $X$.
 \begin{example}\label{exa:45}   Suppose
 $
p_t(x_i,x_i)\le \kappa t^{-d/2}~~{\rm for~}t\ge t_1
{\rm~and~}i\in\{1,2\}.
 $
Suppose $\nu_{x_1}, \nu_{x_2}\ge t_1^{-\kappa}$. Then for each
$\epsilon>0$, there exists a constant $C_0(d,\kappa,\epsilon)>0$
such that
\begin{equation}\label{e:1340we}
p_t(x_1,x_2)\le C_0
t^{-d/2}\exp\left(-\theta\frac{d_\nu(x_1,x_2)^2}{t}
\right)~~{\rm~for~}t\ge t_1^{1+\epsilon}\vee d_\nu(x_1,x_2).
\end{equation}
\end{example}
\pf Let   $f_i(t)=\kappa^{-1}\nu_{x_i}^{-1}t^{d/2}$ for $i\in
\{1,2\}$ and $t\in \bR_+$. Then for each $t\ge t_1$,
$$
\P_{x_i}(X_t=x_i)=\nu_{x_i}p_t(x_i,x_i)\le \frac{1}{f_i(t)}.
$$
Note that   $f_i$ is (1,2)-regular and for $t\ge t_1$,
$$
f_i(t)=\kappa^{-1}\nu_{x_i}^{-1}t^{d/2}\le
\kappa^{-1}t_1^{\kappa}t^{d/2}\le \kappa^{-1} t^{\kappa+d/2}\le A
\exp(2^{-9}t^\epsilon),
$$
where $A$ is some constant which depends only on  $d,\kappa$ and
$\epsilon$. Applying Theorem 5.1 gives
$$
\P_{x_1}(X_{t}=x_2)\le \frac{C_1 (\nu_{x_2}/\nu_{x_1})^{1/2}
}{\sqrt{f_1(\frac{t}{4})f_2(\frac{t}{4})}}\exp\left(-\theta\frac{d_{\nu}(x_1,x_2)^2}{t}\right)~~{\rm~for~}t\ge
t_1^{1+\epsilon}\vee d_\nu(x_1,x_2),
$$
which implies (\ref{e:1340we}) immediately.\qed\\
\\
 {\bf Acknowledgement:} I am  greatly indebted to Martin
Barlow for suggesting the problem and for many stimulating
conversations. I gratefully acknowledge the valuable comments and
discussions with Minzhi Zhao of Zhejiang University. The research
was partially supported by the NSFC grant No. 11001173, and China
Scholarship Council during the author's visit to University of
British Columbia in 2012.

\end{document}